\documentclass[leqno,10pt]{amsart}

\usepackage{latexsym,amssymb,amsthm,amsmath,amscd}

\theoremstyle{plain}
\newtheorem{theorem}{Theorem}[section]

\numberwithin{equation}{section}

 \newtheorem*{Theorem A}{{\bf Theorem A}}
\newtheorem{corollary}{Corollary}[section]

\newtheorem{proposition}{Proposition}

\numberwithin{equation}{section}

\theoremstyle{remark}
\newtheorem{remark}{Remark}[section]

 \numberwithin{equation}{section}
\def\({\left( }
\def\){\right)}
\def\i{{\rm i}}
\def\e{\eqref}

\begin{document}
\title[Geometric classifications of homogeneous production
functions]{Geometric classifications of homogeneous production
functions}

\author[B.-Y. Chen, G.E. V\^{\i}lcu]{Bang-Yen Chen, Gabriel Eduard V\^{\i}lcu}

\date{}
\maketitle

\abstract
In this paper, we completely classify homogeneous production functions with an arbitrary number of inputs whose production hypersurfaces are flat. As an immediate consequence, we obtain a complete classification of homogeneous production functions with two inputs whose
production surfaces are developable.
\\ \\
{\bf Keywords:} Gauss-Kronecker curvature, developable surface, production
function, production hypersurface, return to scale, flat hypersurface.
\\ \\
{\bf 2010 Mathematics Subject Classification:} Primary 91B38; Secondary 53A05.\\
{\bf JEL-codes:} C65, C69.\\
\endabstract

\section{Introduction}

A \emph{developable surface} in Euclidean 3-space $\mathbb E^{3}$ is a smooth surface
$S$ with Gaussian curvature $K\equiv0$
\cite{TOP}. By Gauss' fundamental theorem, a developable surface is
locally isometric to a planar domain and this motivates the name
\emph{developable} \cite{STO2}. Since $K=k_1k_2$, where $k_1$ and
$k_2$ are the principal curvatures, the points of $S$ are either
parabolic or planar points. Such surfaces were classified as the
plane, the cylinder, the cone, or the tangent surface of a space
curve, provided that all points of each surface are either planar
points or all are parabolic points (see e.g. \cite{KLI,SPI,VAI}).
But, despite the fact that there exist many obvious examples of
developable surfaces which contain points of both kinds, it seems
surprising that a long period geometers neglect completely the case
of a surface of curvature zero having both flat and non-flat points.
Indeed, it is an amazing fact in the history of differential
geometry that a strong global result concerning developable surfaces
was announced by A.V. Pogorelov (without proof) only somewhat late
in its development, in 1956, in \cite{POG}: if $S\subset
\mathbb{E}^3$ is a complete surface with zero Gaussian curvature
everywhere, then $S$ must be a generalized cylinder in
$\mathbb{E}^3$, i.e. a surface $S$ that may be described as follows:
there exists a curve $c:\mathbb{R}\rightarrow\mathbb{E}^3$, and a
fixed direction $n$ such that $f(s,t)=c(t)+sn$,
$f:\mathbb{R}^2\rightarrow\mathbb{E}^3$ being a global
parametrization of $S$. Different proofs of this theorem were given
in \cite{DC,HN,MAS,STO}. If the hypothesis of completeness is
omitted, the situation is more complicated, as W.S. Massey remarked
in \cite{MAS} and geometers felt that the possibilities were too
numerous to obtain significant results (see also \cite[p.
355]{WIN}). However, some interesting results can be obtained even
in this case, if we consider developable surfaces described as the
graphs of functions $z=f(x,y)$, where $(x,y)\in D$, $D$ being an
open subset of $\mathbb{R}^2$. It is well known that in this case
the Gaussian curvature of the surface $z=f(x,y)$, where $f$ is
differentiable of class $C^2$, is given by
\[
K=\frac{f_{xx}f_{yy}-f^2_{xy}}{(1+f^2_x+f^2_y)^2},
\]
where the subscripts denote partial derivatives of the function $f$
with respect to the corresponding independent variable. Therefore
the developable surfaces are described as the solutions of the
homogeneous Monge-Amp\`{e}re equation (see \cite[p. 10]{CH})
\begin{equation}\label{e1}
f_{xx}f_{yy}-f^2_{xy}=0.
\end{equation}

The general solutions of equation (\ref{e1}) in parametric form,
under the hypothesis $f_{xx}\neq0$ everywhere, were obtained and
expressed in terms of two functions of one variable $g$ and $h$, by
V.G. Ushakov in \cite{USH}:
\[
\left\{%
\begin{array}{ll}
    x(t,s)=g(t)-sh'(t)  \\
    y(t,s)=s \\
    f(t,s)=tg(t)-\int_0^tg(r)dr+s\left[h(t)-th'(t)\right], \\
\end{array}%
\right.
\]
where $h\in C^2$ and $g\in C^1$ are arbitrary functions such that
$g'\neq0$. His approach is based on the use of a series of changes
of the involved variables deduced by geometric arguments. This
result allows to describe explicitly all developable surfaces
without planar points in $\mathbb{E}^3$. We also remark that the
general solution of the homogeneous Monge-Amp\`{e}re equation in the
higher-dimensional case has been obtained recently by Y. Bozhkov in
\cite{BOZ}, using the same approach.

It is interesting that the motivation for this work comes not only
from classical differential geometry, but also from microeconomics and macroeconomics. It is well-known that the study of the shape and the properties of the production possibility frontier is a subject of great interest in economic analysis. Some conditions under which the production possibility frontier is a surface lying in the edge of a cone, cylinder or plane were derived in \cite{INO1}. Moreover, in \cite{AOT}, the authors investigate the condition for the production surface to contain a flat portion.
On the other hand, in \cite{VV,VGE} it was proved that some
well-known production functions exhibit constant return to scale if
and only if the corresponding hypersurfaces have null Gauss-Kronecker curvature and this
result seems to be true for a large family of classical production
functions found in the economical literature. But the notion of
constant return to scale is equivalent to a homogeneity of degree
$1$ for the production function and therefore a natural question
arises:
{\it  Is it true that a production function has constant return
to scale if and only if the corresponding hypersurface has null Gauss-Kronecker curvature?}
\vskip.1in

The main result of this paper is to prove the following theorem.

\begin{theorem} \label{T:1} Let $f$ be a twice differentiable, $r$-homogeneous, non-constant, real valued function of $n$ variables $(x_{1},\ldots,x_{n})$
on an open domain $D\subset\mathbb{R}^n,\, n\geq 2$.  Then the hypersurface of $\mathbb E^{n+1}$ defined by
$$z=f(x_{1},\ldots,x_{n}),\;\; (x_{1},\ldots,x_{n})\in D,$$ is flat if and only if either $f$ is linearly homogeneous, \emph{i.e.} $r=1$, or
\begin{align}\label{f}f=\left(c_1x_{1}+c_2x_{2}+\cdots+c_{n} x_{n}\right)^r,\;\; r\ne 1,\end{align} for some real constants $c_{1},\ldots,c_{n}$.
 \end{theorem}

It is very important to note that, as we can see in Section 4, the above Theorem is false if we replace the flatness of the hypersurface by vanishing of Gauss-Kronecker curvature. So, somewhat unexpected, we obtain a negative answer to the question naturally raised above. We also remark that Theorem \ref{T:1} extends a result of Brickell concerning homogeneous functions of degree two \cite{BR}.

\section{Some preliminaries on the geometry of hypersurfaces}

For general references on the geometry of hypersurfaces, we refer to \cite{CH1,CH2,REI}.

Let $M$ be a hypersurface of a Euclidean
$(n + 1)$-space $\mathbb E^{n+1}$. The \emph{Gauss map}
$\nu: M \rightarrow S^{n}$ maps $M$ to the unit hypersphere
$S^n$ of $\mathbb E^{n+1}$. This is a
continuous map such that $\nu(p)$ is a unit normal vector $\xi(p)$ of $M$ at $p$. The
Gauss map can always be defined locally, \emph{i.e.} on a small piece of the hypersurface.
It can be defined globally if the hypersurface is orientable.

The differential $d\nu$ of $\nu$ can be used to define a type of extrinsic quantity, known
as the \emph{shape operator}. Since each tangent space $T_pM$ is an inner product space, the
shape operator $S_p$ can be defined as a linear operator on $T_pM$ by
\[
g(S_pv,w)=g(d\nu(v),w)
\]
for $v,w\in T_pM$, where $g$ is the metric tensor on $M$ induced from the Euclidean metric
on $\mathbb E^{n+1}$. The eigenvalues of the shape
operator are called principal curvatures. The determinant of the shape operator
$S_p$, denoted by $K(p)$, is called the \emph{Gauss-Kronecker curvature}. When $n = 2$, the
Gauss-Kronecker curvature is simply called the \emph{Gauss curvature}, which is
intrinsic due to Gauss' theorema egregium.

The Riemann curvature tensor $R$ of $M$ is given in term of the Levi-Civita connection $\nabla$ of $g$ by the following formula:
\[
R(u,v)w=\nabla_u\nabla_vw-\nabla_v\nabla_uw-\nabla_{[u,v]}w.
\]

The curvature tensor measures non-commutativity of the covariant derivative, and
as such is the integrability obstruction for the existence of an isometry with Euclidean
space. In this context, a Riemannian manifold is called \emph{flat} if its Riemann curvature
tensor vanishes identically.

The following basic result is well-known.

\begin{proposition}
For the production hypersurface of $\mathbb E^{n+1}$ defined by
\[
L(x_1,...,x_n)=(x_1,...,x_n,f(x_1,...,x_n)),
\]
we have:\\
(i) The Gauss-Kronecker curvature $K$ is given by
 \begin{equation}\label{p.1}
K=\frac{{ \rm det}(f_{x_ix_j})}{w^{n+2}}
\end{equation}
with $w=\sqrt{1+\sum_{i=1}^{n} f_{i}^{2}}.$\\
(ii) The Riemann curvature tensor $R$ and the metric tensor $g$ satisfy
 \begin{equation}\label{p.2} g\(R\(\frac{\partial}{\partial x_{i}},\frac{\partial}{\partial x_{j}}\)\frac{\partial}{\partial x_{k}},\frac{\partial}{\partial x_{\ell}}\) =\frac{f_{x_ix_{\ell}}f_{x_jx_k}- f_{x_ix_k}f_{x_jx_{\ell}}}{w^{4}}.
\end{equation}
\end{proposition}

\section{Proof of Theorem \ref{T:1}}

Let us assume that $n=2$. We first prove the left-to-right implication. From Euler's
Homogeneous Function Theorem we have
\begin{equation}\label{e2}
xf_{x}+yf_{y}=rf.
\end{equation}

By derivation with respect to $x$ in (\ref{e2}) we obtain
\begin{equation}\label{e3}
(r-1)f_x=xf_{xx}+yf_{yx}
\end{equation}
and similarly, by derivation with respect to $y$ in (\ref{e2}), we
deduce
\begin{equation}\label{e4}
(r-1)f_y=xf_{xy}+yf_{yy}.
\end{equation}

Using now (\ref{e2}), (\ref{e3}) and (\ref{e4}) in the homogeneous
Monge-Amp\`{e}re equation (\ref{e1}) we derive that
\begin{equation}\label{e6}
(r-1)\left[(r-1)f_{x}f_{y}-rff_{xy}\right]=0.
\end{equation}

So, we deduce that $r=1$, i.e. $f$ is homogeneous of degree $1$, or
$f$ satisfy the partial differential equation:
\begin{equation}\label{e7}
(r-1)f_{x}f_{y}-rff_{xy}=0.
\end{equation}

In order to finish the proof of the left-to-right implication, we must
solve (\ref{e7}). But $f$ being homogeneous of degree $r$, it
follows that it can be written in the form:
\[
f(x,y)=y^rh(u)
\]
or
\[
f(x,y)=x^rh(u),
\]
where $u=\frac{x}{y}$ (with $y\neq0$), respectively $u=\frac{y}{x}$
(with $x\neq0$), and $h$ is a real valued function of $u$, of class
$C^2$ on its domain of definition. We can suppose, without loss of
generality, that the first situation occurs, so $f(x,y)=y^rh(u)$,
with $u=\frac{x}{y}$. Then we obtain
\begin{equation}\label{e8}
f_x=y^{r-1}h'(u),
\end{equation}
\begin{equation}\label{e9}
f_y=y^{r-1}[rh(u)-uh'(u)]
\end{equation}
and
\begin{equation}\label{e10}
f_{xy}=y^{r-2}[(r-1)h'(u)-uh''(u)].
\end{equation}

Using (\ref{e8}), (\ref{e9}) and (\ref{e10}) in (\ref{e7}), we
deduce that $h$ must satisfy the differential equation:
\begin{equation}\label{e11}
(r-1)(h'(u))^2-rh(u)h''(u)=0.
\end{equation}

Making now the substitution
\begin{equation}\label{e12}
w(u)=\frac{uh'(u)}{h(u)}
\end{equation}
we derive
\begin{equation}\label{e13}
h'(u)=\frac{w(u)h(u)}{u}
\end{equation}
and we obtain
\begin{equation}\label{e14}
h''(u)=\frac{h(u)}{u}\left(w'(u)-\frac{w(u)}{u}+\frac{w^2(u)}{u}\right).
\end{equation}

Using now (\ref{e13}) and (\ref{e14}) in (\ref{e11}), we deduce that
$w$ must satisfy the following differential equation:
\begin{equation}\label{e15}
uw'(u)-w(u)+\frac{1}{r}w^2(u)=0.
\end{equation}

But (\ref{e15}) is a differential equation of Bernoulli type, having
solution
\begin{equation}\label{e16}
w(u)=\frac{ru}{u+rc},
\end{equation}
where $c$ is a real constant such that $w$ is well defined and of
class $C^2$ on its domain of definition. From (\ref{e12}) and
(\ref{e16}) we derive that $h$ is given by
\[
h(u)=\left(c_1u+c_2\right)^r,
\]
where $c_1$ and $c_2$ are arbitrary real constants. So, we obtain
that $f$ is defined by
\[f(x,y)=\left(c_1x+c_2y\right)^r\]
and the direct implication is proved.

Next, we show that the right-to-left implication also holds for $n=2$.
Indeed, if $r=1$ then it follows from Euler's Homogeneous Function
Theorem that
\begin{equation}\label{e17}
xf_{xx}=-yf_{yx}
\end{equation}
and
\begin{equation}\label{e18}
yf_{yy}=-xf_{xy}.
\end{equation}
But $f$ being differentiable of class $C^2$, we have the equality of
mixed partial derivatives. Hence from (\ref{e17}) and (\ref{e18}),
we deduce that $f$ satisfy homogeneous Monge-Amp\`{e}re equation.
Therefore the surface given by $z=f(x,y)$ is developable.

On the other hand, if $f(x,y)=\left(c_1x+c_2y\right)^r$, then a
direct computation shows that the Gaussian curvature $K$ vanishes
everywhere. Hence the surface described as the graph of function
$z=f(x,y)$ is also developable.
This proves the theorem for $n=2$.

Now, we  assume that $n>2$.
Let $f$ be a twice differentiable real valued function of $n$ variables $(x_{1},\ldots,x_{n})$ defined
on an open domain $D\subset\mathbb{R}^n$, which is  homogeneous of degree $r$. Then we have
   \begin{align}\label{1.1} f(tx_{1},\ldots,tx_{n}) = t^{r}f(x_{1},\ldots,x_{n})\end{align}
   for an $t\in \mathbb{R}$ for which \e{1.1} is defined.
   Since $f$ is  $r$-homogeneous, the Euler Homogeneous Function Theorem implies that
      \begin{align}\label{1.2} x_{1}f_{x_{1}}+x_{2}f_{x_{2}}+\cdots+x_{n} f_{x_{n}}=r f.\end{align}
After taking the partial derivatives of \e{1.2} with respect to $x_{1},\ldots,x_{n}$, respectively, we obtain
  \begin{equation} \begin{aligned}\label{1.3} &x_{1}f_{x_{1}x_{1}}+x_{2}f_{x_{1}x_{2}}+\cdots+x_{n} f_{x_{1}x_{n}}=(r-1) f_{x_{1}},
  \\& x_{1}f_{x_{1}x_{2}}+x_{2}f_{x_{2}x_{2}}+\cdots+x_{n} f_{x_{2}x_{n}}=(r-1) f_{x_{2}},
  \\& \hskip1in \vdots
  \\ &x_{1}f_{x_{1}x_{n}}+x_{2}f_{x_{2}x_{n}}+\cdots+x_{n} f_{x_{n}x_{n}}=(r-1) f_{x_{n}}.\end{aligned}  \end{equation}
  Let us assume that the hypersurface $z=f(x_{1},\ldots,x_{n})$ is flat. Then the hypersurface has vanishing Riemann curvature tensor $R$.
 Because the curvature tensor of the hypersurface satisfies (\ref{p.2}), then we derive from the flatness  that
     \begin{equation} \begin{aligned}\label{1.7} &f_{x_{i}x_{\ell}}f_{x_{j}x_{k}}- f_{x_{i}x_{k}}f_{x_{j}x_{\ell}}=0,\;\; 1\leq i,j,k,\ell\leq n.\end{aligned}  \end{equation}
Therefore, after applying the first two equations of \e{1.3} and by using \e{1.7}, we obtain
  $$(r-1)f_{x_{1}}f_{x_{2}x_{k}}=(r-1)f_{x_{2}}f_{x_{1}x_{k}},\;\; k=1,\ldots,n,$$
which implies that either $r=1$ or $f_{x_{1}}f_{x_{2}x_{k}}=f_{x_{2}}f_{x_{1}x_{k}}$.

  Similarly, by applying the same argument we may conclude that either $r\ne 1$ or $f_{x_{i}}f_{x_{j}x_{k}}=f_{x_{j}}f_{x_{i}x_{k}}$ holds for other $i,j$. Hence we obtain either $r=1$ or
 \begin{align}\label{1.8} f_{x_{i}}f_{x_{j}x_{k}}=f_{x_{j}}f_{x_{i}x_{k}},\;\; 1\leq i\ne j\leq n.\end{align}
  If $r=1$, then $f$ is linearly homogeneous. Therefore from now on we may assume that $r\ne 1$.

\vskip.1in
 {\sc Case} (i): $f_{x_{1}},\ldots,f_{x_{n}}\ne 0$. Now, it follows from \e{1.8} that
 \begin{align}\label{1.9} \frac{f_{x_{1}x_{k}}}{f_{x_{1}}}=\ldots=\frac{f_{x_{n}x_{k}}}{f_{x_{n}}},\;\; k=1\ldots,n.\end{align}
After solving system \e{1.9} for $f_{x_{1}},\ldots,f_{x_{n}}$, we get
 \begin{align}\label{1.10} c_{j}f_{x_{i}}=c_{i}f_{x_{j}},\;\; 1\leq  i\ne j\leq n,\end{align}
for some nonzero real constants $c_{1},\ldots,c_{n}$. Therefore after solving system \e{1.10} we obtain $$f(x_{1},\ldots,x_{n})=F(c_{1}x_{1}+c_{2}x_{2}+\cdots+c_{n}x_{n})$$
for some real-valued function $F$. Since $f$ is assumed to be $r$-homogeneous, we conclude that $f$ is of the form \e{f}.

\vskip.1in
 {\sc Case} (ii): {\it $f_{x_{1}}=\ldots=f_{x_{s}}= 0$ and $f_{x_{s+1}},\ldots,f_{x_{n}}\ne 0$ for some $1\leq s\leq n-2$}. In this case, we derive from \e{1.8} that \begin{align}\label{2.9} \frac{f_{x_{i}x_{k}}}{f_{x_{i}}}=\frac{f_{x_{j}x_{k}}}{f_{x_{j}}},\;\; i,j\in \{s+1,\ldots, n\},\;\; k\in \{1,\ldots,n\}.\end{align}
  Thus, we may apply the same argument as in Case (i) to conclude that
$$f(x_{1},\ldots,x_{n})=(c_{s+1}x_{s+1}+c_{2}x_{2}+\cdots+c_{n}x_{n})^{r},$$
  which is a special case of \e{f}.

\vskip.1in
 {\sc Case} (iii): {\it $f_{x_{1}}=\ldots=f_{x_{n-1}}= 0$ and $f_{x_{n}}\ne 0$}.  In this case, we obtain
 \begin{align}\notag f(x_{1},\ldots,x_{n})=F(x_{n})\end{align}
  for some function $F$. Since $f$ is  $r$-homogeneous, $F$ is the product of a nonzero real constant $c$ times the power function $F(u)=u^{r}$. Consequently, we also obtain a special case of \e{f}.

  The converse is easy to verify.  \qed

An immediate consequence of Theorem \ref{T:1} is the following.

\begin{corollary}\label{C:1}
Let $f$ be a $r$-homogeneous, differentiable of class $C^2$, non-constant, real valued function of two variables $(x,y)$ defined
on an open domain $D\subset\mathbb{R}^2$. Then the surface defined by
$z=f(x,y)$, $(x,y)\in D$, is developable if and only if either $f$ is linearly homogeneous, \emph{i.e.} $r=1$, or
$f(x,y)=\left(c_1x+c_2y\right)^r$, $r\ne 1$, where $c_1$ and $c_2$ are real
constants.
\end{corollary}

\begin{remark} If the domain $D$ has boundary or the domain is a cone, then the same result holds by applying the continuity condition of the production function $f$.
\end{remark}

\section{An application to the theory of production functions}

Almost all economic theories presuppose a production function, either on the firm level or the aggregate level. Therefore, the production functions are a key concept
both in microeconomics and macroeconomics. Roughly speaking, they
are a mathematical formalization of the relationship between the
output of a firm, an industry, or an entire economy, and the inputs
that have been used in obtaining it. Mathematically, a production
function is a map $f$ of class $C^\infty$,
$f:\mathbb{R}_+^n\rightarrow\mathbb{R}_+$, $f=f(x_1,x_2,...,x_n)$,
where $f$ is the quantity of output, $n$ is the number of the inputs
and $x_1,x_2,...,x_n$ are the factor inputs (such as labor, capital,
land, raw materials etc.). We remark that a production function $f$
can be identified with the graph of $f$, \emph{i.e.} the nonparametric hypersurface of the $(n+1)$-dimensional
Euclidean space $\mathbb{E}^{n+1}$, defined by
\[L(x_1,...,x_n)=(x_1,...,x_n,f(x_1,...,x_n))\]
and called the \emph{production hypersurface} of $f$. It is clear that in
the particular case of two inputs, we have a surface. In order for
these functions to model the economic reality, they are required to
have certain properties (see e.g. \cite{FF,SHE,THO}). We recall now
some of them with appropriate economic interpretations:
\begin{enumerate}
\item[$1.$] $f$ vanishes in the absence of an input; this means that the factor
inputs are necessary.
\item[$2.$] $f_{x_i}>0$, for all
$i\in\{1,...,n\}$, which indicates that the production function is
strictly increasing with respect to any factor of production.
\item[$3.$] $f_{x_ix_i}<0$, for all
$i\in\{1,...,n\}$, which signifies that the production has
decreasing efficiency with respect to any factor of production.
\item[$4.$] $f(x+y)\geq f(x)+f(y)$, $\forall x,y\in\mathbb{R}_+^n$,
which means that the production has non-decreasing global
efficiency.
\item[$5.$] $f$ is a homogeneous function, i.e. there exists a real number $r$ such that $f(\lambda\cdot x)=\lambda^r f(x)$ for all $x\in\mathbb{R}_+^n$
and $\lambda\in\mathbb{R}_+$, which signifies that if the inputs are
multiplied by same factor, then the output is multiplied by some
power of this factor. If $r=1$ then the function is said to have a
constant return to scale, if $r>1$ then we have an increased return
to scale and if $r<1$ then we say that the function has a decreased
return to scale.
\end{enumerate}

In 1928, C.W. Cobb and P.H. Douglas \cite{CD} introduced a famous
two-factor production function, nowadays called Cobb-Douglas
production function, in order to describe the distribution of the
national income by help of production functions. This function was
further generalized by K.J. Arrow, H.B. Chenery, B.S. Minhas and
R.M. Solow \cite{ACMS}, they introducing the so-called Constant
Elasticity of Substitution production function. In 1963 it was
generalized to the $n$-factor case by H. Uzawa \cite{UZ} and D.
McFadden \cite{MCF}. This function, usually called generalized CES
production function, Armington aggregator or ACMS function, is
defined by $f:\mathbb{R}_+^n\rightarrow\mathbb{R}_+$,
\begin{equation}\label{9}
f(x_1,...,x_n)=A\left(\sum_{i=1}^{n}c_i
x_i^\rho\right)^{\frac{\gamma}{\rho}},
\end{equation}
where $A>0$, $\rho< 1$, $\rho\neq 0$, $\gamma>0$ and $c_i> 0$, for
all $i\in\{1,...,n\}$.

We note that we can derive some well-known production function from
the generalized CES production function, as special cases:
\begin{enumerate}
\item[$i.$] If we take $c_1>0,...,c_n>0$ such that $\sum_{i=1}^{n}c_i=1$, $\gamma=1$ and $\rho \rightarrow 0$ in (\ref{9}), then  we obtain the Cobb-Douglas production function,
also known as the \emph{imperfect complements} production function:
\[f(x_1,...,x_n)=A\cdot
\prod_{i=1}^n x_i^{c_{i}}.\] If a production function  has the above
expression, but $\sum_{i=1}^{n}c_i$ is not necessary 1, then $f$ is
called a generalized Cobb-Douglas production function.
\item[$ii.$] If we take $\rho \rightarrow 1$ in (\ref{9}), then we obtain the \emph{multinomial} production function:
\[f(x_1,...,x_n)=A\left(\displaystyle\sum_{i=1}^nc_i x_i\right)^\gamma.\]
When $n=2$, the multinomial production function is called \emph{binomial}.
In particular, if $\gamma=1$, then we obtain the \emph{linear} production
function, also called the\emph{ perfect substitutes} production function
(see also \cite{BHA} for other remarkable classes of production
functions).
\end{enumerate}

In \cite{VV,VGE} it was proved that generalized Cobb-Douglas
production functions and generalized CES production functions have
constant return to scale if and only if the corresponding
hypersurfaces have vanishing Gauss-Kronecker curvature. Therefore, it is natural to ask if
a general result of this kind holds for all production functions.
The answer, in the case of a production function with two inputs,
follows now from Corollary \ref{C:1}.

\begin{theorem}\label{2} A homogeneous production function with two inputs defines a production surface with vanishing Gauss curvature
if and only if either it has constant return to scale or it is a binomial production function.
\end{theorem}

We note that the above Theorem completely classifies homogeneous production functions with two inputs whose
production surfaces are developable. On the other hand, from Theorem \ref{T:1} we can also obtain a complete  classification of homogeneous production functions with an arbitrary number of inputs whose production hypersurfaces are flat, as follows.

\begin{theorem}\label{3} A homogeneous production function with an arbitrary number of inputs defines a flat hypersurface if and only if either it has constant return to scale or it is a multinomial production function.
\end{theorem}

\begin{remark}\label{R:1} Theorem \ref{3} is false if the flatness of the hypersurface were replaced by vanishing of Gauss-Kronecker curvature. This can be seen from the following example.
Consider the $r$-homogeneous production function
\begin{align}\label{1.12}f(x,y,z)=(x+y+\sqrt{yz})^{r},\;\; (x,y,z)\in {\mathbb R}_{+}^{3},\;\; r>1.\end{align}
Then the hypersurface in $\mathbb E^{4}$ defined by \e{1.12} is non-flat, but it has vanishing Gauss-Kronecker curvature. We note that recently, in \cite{CH4}, the first author has completely classified quasi-sum production functions whose production hypersurfaces have vanishing Gauss-Kronecker curvature. Other classification results concerning production functions were proved recently in \cite{CH55,CH66,VV2}.
  \end{remark}

\begin{remark}The first author would like to point out that the proof of Theorem 3.1 of \cite{CH3} contains an error. Consequently, Theorem 3.1 of \cite{CH3} shall be replaced by Theorem \ref{3} of this article (see also \cite{CH5}).
\end{remark}

\section{Conclusions}

In this paper we obtain two classification results concerning homogeneous production functions. The first result, given in Theorem \ref{2}, is rather unexpected because there are dozens of classes of homogeneous production functions with two inputs used in
economy (see e.g. \cite{BHA}), but only one can define a production surface having null Gauss curvature, without
exhibit constant return to scale. It is very interesting that this statement does not remain valid for any dimension, as follows from Remark \ref{R:1}. However, the result can be generalized for an arbitrary number of inputs considering the homogeneous production functions that define flat production hypersurfaces, as we can see in Theorem \ref{3}.

\section*{Acknowledgement} The second author was supported by
CNCS-–UEFISCDI, project number PN--II--ID--PCE--2011--3--0118.

\noindent    Bang-Yen CHEN\\
    Department of Mathematics,\\
    Michigan State University,\\
    East Lansing, Michigan
    48824--1027, USA\\
    E-mail address: bychen@math.msu.edu\\ \\

\noindent Gabriel Eduard V\^{I}LCU$^{1,2}$ \\
      $^1$University of Bucharest,\\
      Research Center in Geometry, Topology and Algebra,\\
      Str. Academiei, Nr. 14, Sector 1,\\
      Bucure\c sti 70109, ROMANIA\\
      E-mail address: gvilcu@gta.math.unibuc.ro\\
      $^2$Petroleum-Gas University of Ploie\c sti,\\
      Department of Mathematical Modelling, Economic Analysis and Statistics,\\
      Bd. Bucure\c sti, Nr. 39,\\
      Ploie\c sti 100680, ROMANIA\\
      E-mail address: gvilcu@upg-ploiesti.ro

\end{document}